\documentclass[11pt]{amsart}
\usepackage{amsmath,amssymb,amsthm}
\usepackage{xcolor}
\usepackage{booktabs}
\usepackage{needspace}
\usepackage{microtype}
\usepackage[colorlinks=true,citecolor=blue!45!black,linkcolor=blue!45!black,urlcolor=blue!45!black]{hyperref}
\usepackage{mathtools}

\newtheorem{theorem}{Theorem}[section]
\newtheorem{lemma}[theorem]{Lemma}
\newtheorem{corollary}[theorem]{Corollary}

\theoremstyle{remark}
\newtheorem{remark}[theorem]{Remark}
\newcommand{\qbinom}[2]{\genfrac{[}{]}{0pt}{}{#1}{#2}_{\!q}}
\newcommand{\Fq}{\mathbb F_q}
\newcommand{\rk}{\operatorname{rk}}
\newcommand{\rad}{\operatorname{rad}}

\title[Cross-intersection in semilattices]{An incidence approach to cross-intersection in semilattices}
\author{Sho Suda}
\address{Department of Mathematics, National Defense Academy of Japan, 1-10-20 Hashirimizu, Yokosuka, Kanagawa 239-8686, Japan}

\email{ssuda@nda.ac.jp}
\date{\today}
\hypersetup{pdftitle={An incidence approach to cross-intersection in semilattices},pdfauthor={Sho Suda}}

\begin{document}

\begin{abstract}
We present an incidence-based method for product bounds on cross-intersecting
families in finite graded meet-semilattices. The two families may lie in
different rank or type fibres and may have different design strengths. The
criterion separates global incidence parameters from two directional local
counts and yields the sharp product bound, with equality only for full stars
with a common centre.

Applications include mixed subspace types relative to a fixed subspace, a
type-restricted generator family in $Q^-(7,q)$, and mixed form ranks in
symplectic spaces. The same short criterion gives several sharp results not
obtained directly from the specific theorems and general criteria compared
here. In the singular-space case, it also applies in an infinite parameter
range beyond the natural star-ratio and weak-spreadness tests.
\end{abstract}
\maketitle
\raggedbottom

\section{Introduction}
In a cross-intersection problem, every pair $(A,B)\in\mathcal A\times\mathcal B$
is required to satisfy a prescribed intersection condition, while no
restriction is imposed on pairs within $\mathcal A$ or within $\mathcal B$.
The associated product problem asks for the maximum of $|\mathcal A||\mathcal B|$.
For families of sets, the cross-intersection requirement is typically
$|A\cap B|\ge s$, while for families of subspaces it takes the form
$\dim(A\cap B)\ge s$. In both settings, the obvious candidates for
extremality are \,full stars\, with a shared centre. Beginning with the
product theorems of Pyber \cite{Pyber1986} and Matsumoto--Tokushige
\cite{MT1989}, sharp bounds and stability results have been obtained for
set systems \cite{BaoJi2026,Borg2016,FLST2014} and for vector spaces
\cite{CLLW2023,ST2014,WenLv2026,YLW2025}, including the case where the two
families live in different dimensions.

Our aim is to provide a single local incidence argument for the product bounds
in the setting of finite graded meet-semilattices. Set systems and subspace
families correspond to inclusion posets, while families of partial maps fit
via restriction, so these examples can be treated uniformly.

A class $\mathcal P$ of lower-rank objects controls the notion
of intersection and, at the same time, specifies the allowable centres of
stars. The two ambient families may lie in different rank or type fibres,
and they may be designs of different strengths.

For each side $a\in\{1,2\}$, let $v_a$ denote the number of admissible centres
contained in a member, let $L_a$ be the size of a full star, and let $d_a$
be an upper bound on the number of star members that can remain compatible
with an opposite member that avoids the centre. Theorem~\ref{thm:criterion}
shows that
$ v_ad_a<L_a$ $(a=1,2)$ implies  
\[|\mathcal A||\mathcal B|\le L_1L_2,
\]
and that equality occurs if and only if $\mathcal A$ and $\mathcal B$ are the
two full stars with a common centre. All applications share this extremal
argument; the underlying geometry enters only through the two directional
estimates.

This work continues the local-counting line of research around the
Erd\H{o}s--Ko--Rado (EKR) theorem for designs. Delsarte \cite{Delsarte1976}
introduced designs in regular semilattices, Rands \cite{Rands1982} proved an
EKR theorem for subfamilies of block designs, and Fu \cite{Fu1999}---using
regular quantum matroids in Terwilliger's framework \cite{Terwilliger1996}---
covered Grassmann and bilinear-forms designs. Suda \cite{Suda2012} extended
this approach to a wider class of semilattices, and Guo, Ma and Wang
\cite{GMW2013} developed the regularized and partition-regularized
frameworks.

Our criterion separates the cross-product argument from the regularity
assumptions used to compute its parameters. When centres are unrestricted,
a rank identity for joins gives sharper local estimates and hence stronger
design bounds. The semilattice and design results cited above concern
one-family intersection problems and regularity conditions; they do not
give the mixed-fibre cross-product criterion used here.

Other approaches address related intersection problems.
Algebraic and spectral methods provide a powerful framework for
EKR-type problems in highly symmetric settings; see, for example, Godsil and Meagher~\cite{GodsilMeagher2016}. Qian and Ray-Chaudhuri \cite{QR1997} apply polynomial and incidence-matrix
techniques to a single family under restrictions on meet ranks, including
modular constraints. Kiermaier and Klawuhn \cite{KK2026} derive
intersection-number identities and bounds for designs in a weakened
regular-semilattice setting. For set systems, general cross-product
criteria include Borg's star-ratio theorem
\cite[Theorem~1.2]{Borg2017} and the weak-spreadness criterion of
Wen and Lv \cite[Theorem~1.3(ii)]{WLUnified2026}.

In contrast, the present method uses a two-sided local count that depends on the geometry of the \emph{opposite} family, rather than only on the parameters of each ambient family separately.

Three mixed-type applications illustrate this distinction. They also show that the criterion gives several sharp product results that do not follow directly from the general or geometric results compared below.

\begin{itemize}
\item \emph{Singular spaces.} 
Theorem~\ref{thm:singular} pairs $k$-subspaces disjoint from a fixed
$\ell$-subspace $W$ with $(c+1)$-subspaces meeting $W$ in dimension $c$.
In one of the two local counts, projection to $V/W$ leaves at most one
possible intersection point. When $\ell=2c$, this gives an infinite
parameter range in which the natural outside-point encoding does not
satisfy the hypotheses of either general criterion above; see
Section~\ref{subsec:borg}.

\item \emph{Dual polar spaces.} 
Corollary~\ref{cor:dual-polar-type} considers generators of $Q^-(7,q)$
that meet a fixed hyperbolic six-space in a point, with permitted centres
given by the totally singular planes disjoint from that six-space.
The direct local count is one, yielding a sharp product bound that does
not follow directly from the generator results compared below. 

\item \emph{Symplectic spaces.} 
Theorem~\ref{thm:mixed-symplectic} pairs four-dimensional subspaces
of form ranks $4$ and $2$ in a symplectic space, requiring each
cross-intersection to contain a nondegenerate plane. The local counts
depend on the position of the intersection point relative to the radical.
\end{itemize}

The related geometries have been studied in several settings; see
\cite{OLW2014,ZLL2016} for singular spaces,
\cite{Ihringer2015,IhringerMetsch2018} for generator intersection problems
in classical polar spaces and dual polar graphs, and
\cite{YWAffine,YWPolar} for cross-intersection problems in affine and
totally isotropic symplectic settings.  
The singular-space results do not give the mixed cross-product theorem, 
the dual-polar results do not cover the type-restricted generator family, and
the affine and symplectic-polar results do not cover the mixed form-rank pair in Theorem~\ref{thm:mixed-symplectic}. More detailed comparisons are given in Section~\ref{sec:mixed-types}.
We do not claim that these parameter conditions are optimal.

Section~\ref{sec:criterion} proves the criterion, and
Section~\ref{sec:regularized} develops semilattices and designs.
The mixed-type applications appear in Section~\ref{sec:mixed-types}, and 
Section~\ref{sec:applications} summarizes standard full-fibre consequences, including their interpretation for attenuated spaces. 

Throughout, $q$ is a prime power, and points and planes have vector-space
dimensions one and two, respectively. We use the Gaussian coefficient
\[
 \qbinom ab=\prod_{j=0}^{b-1}\frac{q^{a-j}-1}{q^{b-j}-1}
 \qquad(0\le b\le a).
\]

\section{The local-incidence criterion}\label{sec:criterion}
Let $(X,\preceq)$ be a finite graded meet-semilattice with rank function
$\rk$ and maximum rank $m$, and write $X_r=\{x\in X:\rk(x)=r\}$.
Choose nonempty classes
\[
 \mathcal F_a\subseteq X_{N_a}\quad(a=1,2),\qquad
 \mathcal P\subseteq X_s,\qquad 1\le s\le N_a\le m.
\]
The classes $\mathcal F_a$ may be full rank fibres or prescribed
type classes, and $\mathcal P$ may be taken as $X_s$ in general or $X_{s,i}$ in the partition-regularized framework in Section~\ref{sec:regularized}.  
For $Y\subseteq X$ and $x\in X$, set
$Y_x=\{z\in Y:x\preceq z\}$; if $x\in\mathcal P$, this is a
\emph{$\mathcal P$-star}. Families $Z_a\subseteq\mathcal F_a$ are
\emph{cross $\mathcal P$-intersecting} if every pair
$(z_1,z_2)\in Z_1\times Z_2$ has a common lower bound in $\mathcal P$. 
Throughout, all families in cross-intersection problems are assumed to be nonempty. 
For $\mathcal P=X_s$, this is equivalent to
$\rk(z_1\wedge z_2)\ge s$, and we say \emph{cross $s$-intersecting}.

Assume that each member of $\mathcal F_a$ contains exactly $v_a$
elements of $\mathcal P$, where $v_a>0$. 
Let $Y_a\subseteq\mathcal F_a$ be nonempty for $a=1,2$,
and suppose that there exist positive integers $L_1,L_2$
such that
\begin{equation}\label{eq:balanced}
 |(Y_a)_x|=L_a\qquad(x\in\mathcal P,\ a=1,2).
\end{equation}
Thus each $Y_a$ is an incidence design relative to $\mathcal P$.
For $a\in\{1,2\}$, write $\bar a=3-a$ for the other index. 
For $a\in\{1,2\}$, let $d_a\ge0$ be such that,
for every $x\in\mathcal P$ and $y\in Y_{\bar a}$
with $x\npreceq y$,
\begin{equation}\label{eq:local-count}
 \left|\left\{z\in(Y_a)_x:
       \exists u\in\mathcal P\text{ with }u\preceq z\wedge y
 \right\}\right|\le d_a.
\end{equation}


\begin{theorem}[Mixed-fibre cross-intersection criterion]\label{thm:criterion}
Suppose that $v_ad_a<L_a$ for $a=1,2$.
Then any cross $\mathcal P$-intersecting families
$Z_a\subseteq Y_a$ satisfy
\[
 |Z_1||Z_2|\le L_1L_2.
\]
Equality holds if and only if there exists $x\in\mathcal P$
such that $Z_a=(Y_a)_x$ for $a=1,2$.
\end{theorem}

\begin{proof}

Put $n_a(x)=|(Z_a)_x|$. Counting the sets
\[
 \{(x,z_1,z_2)\in\mathcal P\times Z_1\times Z_2:
   x\preceq z_1\wedge z_2\}, \quad \text{and}\quad
 \{(x,z)\in\mathcal P\times Z_a:x\preceq z\}
\]
in two ways yields
\begin{equation}\label{eq:count}
 |Z_1||Z_2|\le\sum_{x\in\mathcal P}n_1(x)n_2(x), \quad \text{and}\quad
 \sum_{x\in\mathcal P}n_a(x)=v_a|Z_a|.
\end{equation}


\emph{Case 1.} Neither family is contained in a $\mathcal P$-star.
Then, for each $x\in\mathcal P$ and $a\in\{1,2\}$, there exists
$y\in Z_{\bar a}$ with $x\npreceq y$.
Hence $n_a(x)\le d_a$ for all $x\in\mathcal P$ and $a=1,2$.
By \eqref{eq:count},
\[
 |Z_1||Z_2|
 \le \sum_{x\in\mathcal P}n_1(x)n_2(x)
 \le d_1\sum_{x\in\mathcal P}n_2(x)
 =v_2d_1|Z_2|.
\]
Since $Z_2$ is nonempty, we obtain $|Z_1|\le v_2d_1$.
Interchanging the roles of the two families gives
$|Z_2|\le v_1d_2$. Therefore
\[
 |Z_1||Z_2|
 \le (v_1d_1)(v_2d_2)
 <L_1L_2.
\]

\emph{Case 2.} At least one family is contained in a $\mathcal P$-star.
By symmetry, we may assume that
$Z_1\subseteq(Y_1)_x$ for some $x\in\mathcal P$.
If $Z_2\subseteq(Y_2)_x$ as well, the bound and equality case are immediate.
If not, a member of $Z_2$ missing $x$ gives $|Z_1|\le d_1$.
Fix $z\in Z_1$. Every member of $Z_2$ contains one of the $v_1$ centres
below $z$. Hence $|Z_2|\le v_1L_2$ and
$|Z_1||Z_2|\le v_1d_1L_2<L_1L_2$.
\end{proof}

We next give a general upper bound for $d_a$.
In a finite meet-semilattice, two elements with a common upper bound have a join.

\begin{lemma}\label{lem:extension}
Suppose that $s<m$, and set
\[
 M_a=\max_{w\in X_{s+1}}|(Y_a)_w|.
\]
Then one may take $d_a=v_{\bar a}M_a$
in Theorem~\ref{thm:criterion}.
\end{lemma}
\begin{proof}
Fix $x\in\mathcal P$ and $y\in\mathcal F_{\bar a}$ with $x\npreceq y$.
For each $z$ in the set counted in \eqref{eq:local-count}, choose $u\in\mathcal P$ with $u\preceq z\wedge y$. 
Then $u\ne x$, and $x\vee u$ exists and has rank at least $s+1$.
For each such pair $(x,u)$, fix an element $w(x,u)\in X_{s+1}$ satisfying $x\preceq w(x,u)\preceq x\vee u$. 
Since $x\vee u\preceq z$, we have $w(x,u)\preceq z$. 
There are at most $v_{\bar a}$ choices for $u$,
and each $w(x,u)$ is contained in at most $M_a$
members of $Y_a$. 
Therefore the counted set has size at most $v_{\bar a}M_a$. 
\end{proof}

\begin{corollary}\label{cor:balanced}
Under the preceding incidence assumptions, suppose that $s<m$ and $v_1v_2M_a<L_a$ for $a=1,2$. 
Then the bound and equality characterization of Theorem~\ref{thm:criterion} hold.
\end{corollary}


\section{Regularized semilattices and designs}\label{sec:regularized}
We relate our criterion to the design methods of
\cite{Fu1999,Rands1982,Suda2012} and the partition-regularized
framework of \cite{GMW2013}.
We state the incidence regularity assumptions explicitly,
since the mixed-rank setting requires counts relative to both fibres.

\subsection{Partitioned fibres}
Suppose that $s<m$ and that each rank fibre is partitioned into types $X_{r,i}$, and set
\[
 \mathcal P=X_{s,i},\qquad \mathcal F_a=X_{N_a,j_a}.
\]
Assume that, for the ranks and types considered below, the quantities
\begin{align*}
 v_a&=|\{x\in X_{s,i}:x\preceq z\}| \qquad (z\in\mathcal F_a),\\
 \theta_a(r,j)&=|\{z\in\mathcal F_a:w\preceq z\}| \qquad (w\in X_{r,j}), 
\end{align*}
are independent of the choices of $z$ and $w$, respectively.
We also assume that $v_a>0$ and $\theta_a(s,i)>0$.
These are the only incidence regularity conditions needed below. 
Suppose that
\begin{equation}\label{eq:prs-condition}
 v_1v_2\max_j\theta_a(s+1,j)<\theta_a(s,i)
 \qquad(a=1,2),
\end{equation}
where the maximum is taken over the nonempty types $X_{s+1,j}$.
Then Corollary~\ref{cor:balanced} gives
\[
 |Z_1||Z_2|\le\theta_1(s,i)\theta_2(s,i)
\]
for all cross $\mathcal P$-intersecting families
$Z_a\subseteq\mathcal F_a$.
Equality holds if and only if there exists $x\in X_{s,i}$ such that $Z_a=(X_{N_a,j_a})_x$ for $a=1,2$. 
For incidence designs $Y_a\subseteq\mathcal F_a$,
replace $\theta_a(s,i)$ and $\max_j\theta_a(s+1,j)$
by $L_a$ and $M_a$, respectively.

\subsection{Unrestricted centres and the join identity}
When $\mathcal P=X_s$, the rank identity for joins gives a sharper local estimate. 
Assume that every pair $x,u$ with a common upper bound satisfies
\begin{equation}\label{eq:rank-join}
 \rk(x\vee u)=\rk(x)+\rk(u)-\rk(x\wedge u).
\end{equation}
Suppose also that, for $k\in\{N_1,N_2\}$ and
$0\le r\le j\le k$, the quantities
\begin{align*}
 \nu(j,k)&=|\{u\in X_j:u\preceq y\}|,
       &&y\in X_k,\\
 \theta(r,k)&=|\{z\in X_k:x\preceq z\}|,
       &&x\in X_r,\\
 \mu(r,j,k)&=|\{u\in X_j:x\preceq u\preceq y\}|,
       &&x\in X_r,\ y\in X_k,\ x\preceq y,
\end{align*}
are independent of $y$, $x$, and the pair $(x,y)$, respectively.
We further assume that $\theta(s,N_a)>0$ for $a=1,2$.

\begin{corollary}\label{cor:srs}
Under these assumptions, let $1\le s<\min\{N_1,N_2\}$,
and set $v_a=\nu(s,N_a)$ and $L_a=\theta(s,N_a)$.
Suppose that
\begin{equation}\label{eq:srs-condition}
 v_a\max_{\max\{0,2s-N_a\}\le r<s}
 \mu(r,s,N_{\bar a})\theta(2s-r,N_a)<L_a
 \qquad(a=1,2). 
\end{equation}
Then any cross $s$-intersecting families
$Z_a\subseteq X_{N_a}$ satisfy
\[
 |Z_1||Z_2|\le L_1L_2.
\]
Equality holds if and only if 
 there exists $x\in X_{s}$ such that $Z_a=(X_{N_a})_x$
for $a=1,2$.
\end{corollary}
\begin{proof}
Fix $x\in X_s$ and $y\in X_{N_{\bar a}}$
with $x\npreceq y$, and put
$r=\rk(x\wedge y)<s$. 
For each $z$ counted in \eqref{eq:local-count},
choose $u\in X_s$ such that
$x\wedge y\preceq u\preceq z\wedge y$.
Then $x\wedge u=x\wedge y$, and \eqref{eq:rank-join} gives
$\rk(x\vee u)=2s-r$.

If $2s-r>N_a$, no such $z$ exists.
Otherwise, there are at most $\mu(r,s,N_{\bar a})$
choices for $u$, and each $x\vee u$ is contained in
at most $\theta(2s-r,N_a)$ members of $X_{N_a}$.
Taking the maximum over the possible values of $r$
gives the required bound on $d_a$.
The conclusion follows from Theorem~\ref{thm:criterion}.
\end{proof}

\subsection{Designs of different strengths}
The same argument yields a version for designs by
truncating the extension rank at the design strength,
as in the one-family case of \cite{Suda2012}.

\begin{corollary}[Designs of different strengths]\label{cor:designs}
Under the same regularity assumptions, for $a=1,2$,
let $Y_a\subseteq X_{N_a}$ be a nonempty $t_a$-design
with $s<t_a\le N_a$.
For $0\le j\le t_a$, let $\lambda_{a,j}$ denote
the number of members of $Y_a$ containing a given
element of $X_j$, where $\lambda_{a,t_a}>0$. 
Suppose that
\begin{equation}\label{eq:design-condition}
 \nu(s,N_a)
 \max_{\max\{0,2s-N_a\}\le r<s}
 \mu(r,s,N_{\bar a})
 \lambda_{a,\min\{2s-r,t_a\}}
 <\lambda_{a,s}
 \qquad(a=1,2).
\end{equation}
Then any cross $s$-intersecting families
$Z_a\subseteq Y_a$ satisfy
\[
 |Z_1||Z_2|\le\lambda_{1,s}\lambda_{2,s}.
\]
Equality holds if and only if there exists $x\in X_s$
such that $Z_a=(Y_a)_x$ for $a=1,2$.
\end{corollary}
\begin{proof}
Fix $x\in X_j$ with $0\le j\le t_a$, and let $c_x$
denote the number of elements of $X_{t_a}$ containing $x$.
Counting triples $x\preceq u\preceq z$ with $u\in X_{t_a}$,
first for $z\in X_{N_a}$ and then for $z\in Y_a$, gives
\[
 c_x\theta(t_a,N_a)
 =\theta(j,N_a)\mu(j,t_a,N_a)
\]
and
\[
 c_x\lambda_{a,t_a}
 =|(Y_a)_x|\mu(j,t_a,N_a).
\]
Eliminating $c_x$ yields
\[
 |(Y_a)_x|
 =\lambda_{a,t_a}
 \frac{\theta(j,N_a)}{\theta(t_a,N_a)}
 =\lambda_{a,j}
 \qquad(j\le t_a).
\]

For the local bound, a join of rank $k=2s-r$
has exactly $\lambda_{a,k}$ extensions in $Y_a$
if $k\le t_a$.
If $k>t_a$, fix a rank-$t_a$ element below the join.
Every extension of the join in $Y_a$ contains this element,
so there are at most $\lambda_{a,t_a}$ such extensions.
The conclusion follows from
Theorem~\ref{thm:criterion}.
\end{proof}

\section{Mixed subspace types}\label{sec:mixed-types}
In the following three applications, direct local counts yield sharper estimates than the general extension bound. We work in the full subspace lattice, but restrict the two families and the permitted centres to specified geometric types. The star sizes are therefore computed within these restricted families.

\subsection{Singular linear spaces}\label{subsec:singular}
Fix an $\ell$-dimensional subspace $W$ of $V=\Fq^{n+\ell}$. A subspace $U$ has
\emph{type $(r,j)$} if $\dim U=r$ and $\dim(U\cap W)=j$; see
\cite{OLW2014,ZLL2016}. Set
\[
 \begin{split}
 \mathcal F_1&=\{A\le V:\dim A=k,\ A\cap W=\{0\}\},\\
 \mathcal F_2&=\{B\le V:\dim B=c+1,\ \dim(B\cap W)=c\}.
 \end{split}
\]

\begin{theorem}\label{thm:singular}
Let $k\ge2$, $c\ge1$, $n\ge k$, and $\ell\ge2c$.
Suppose that
\begin{equation}\label{eq:singular-condition}
 q^\ell\frac{q^{n-1}-1}{q^{k-1}-1}
   >q^c\qbinom k1.
\end{equation}
If families $\mathcal A\subseteq\mathcal F_1$ and
$\mathcal B\subseteq\mathcal F_2$ satisfy
$A\cap B\ne\{0\}$ for every
$A\in\mathcal A$ and $B\in\mathcal B$, then
\begin{equation}\label{eq:singular-bound}
 |\mathcal A||\mathcal B|
 \le q^{\ell(k-1)}\qbinom{n-1}{k-1}\qbinom\ell c.
\end{equation}
Equality holds if and only if there exists a point
$P\nsubseteq W$ such that
$\mathcal A=(\mathcal F_1)_P$ and
$\mathcal B=(\mathcal F_2)_P$.
In particular, \eqref{eq:singular-condition} holds whenever
$n+\ell\ge2k+c$.
\end{theorem}

\begin{proof}
Let $\mathcal P$ consist of the points outside $W$, and set
$Y_a=\mathcal F_a$ for $a=1,2$.
Since every member of $\mathcal F_1$ is disjoint from $W$,
the intersection condition is equivalent to cross
$\mathcal P$-intersection. The incidence numbers and star sizes are
\begin{equation}\label{eq:singular-counts}
 v_1=\qbinom k1, \quad 
 L_1=q^{\ell(k-1)}\qbinom{n-1}{k-1}, \quad v_2=q^c, \quad 
 L_2=\qbinom\ell c.
 \end{equation}
Fix $P\in\mathcal P$.
The subspace $(W+P)/P$ of $V/P$ has dimension $\ell$.
The members of $(\mathcal F_1)_P$ correspond to graphs of linear
maps from $(k-1)$-subspaces of a fixed $(n-1)$-dimensional
complement of $(W+P)/P$.
There are $\qbinom{n-1}{k-1}$ choices for the domain and
$q^{\ell(k-1)}$ maps for each choice, giving $L_1$.
Every member of $(\mathcal F_2)_P$ can be written uniquely as
$P+D$, where $D$ is a $c$-subspace of $W$, giving $L_2$.

For the first local count, fix $B\in\mathcal F_2$ with
$P\nsubseteq B$.
Every $A\in(\mathcal F_1)_P$ with $A\cap B\ne\{0\}$
contains one of the $q^c$ points of $B$ outside $W$.
For each such point $Q$, if $P+Q$ meets $W$ nontrivially,
no member of $\mathcal F_1$ contains $P+Q$.
Otherwise, exactly
$q^{\ell(k-2)}\qbinom{n-2}{k-2}$ members of $\mathcal F_1$
contain $P+Q$.
Hence we may take
\begin{equation}\label{eq:singular-d1}
 d_1=q^c q^{\ell(k-2)}\qbinom{n-2}{k-2}.
\end{equation}

For the second local count, fix $A\in\mathcal F_1$ with
$P\nsubseteq A$, and let $\pi:V\to V/W$ be the quotient map.
Every $B\in(\mathcal F_2)_P$ satisfies $\pi(B)=\pi(P)$.
Since $\pi|_A$ is injective, there is at most one point
$Q\le A$ with $\pi(Q)=\pi(P)$.
If no such $Q$ exists, the local count is zero.
Otherwise, any nonzero intersection $A\cap B$ must equal $Q$.
Since $Q\ne P$ and $\pi(Q)=\pi(P)$,
the subspace $R=(P+Q)\cap W$ is a point.
Writing $B=P+D$, where $D$ is a $c$-subspace of $W$,
we have $Q\le B$ if and only if $R\le D$.
Thus we may take
\begin{equation}\label{eq:singular-d2}
 d_2=\qbinom{\ell-1}{c-1}.
\end{equation}

The inequality $v_1d_1<L_1$ is equivalent to
\eqref{eq:singular-condition}.
For the other inequality,
\[
 \frac{L_2}{d_2}
 =\frac{q^\ell-1}{q^c-1}
 \ge q^c+1>v_2,
\]
since $\ell\ge2c$.
The bound and equality characterization therefore follow from
Theorem~\ref{thm:criterion}.

Finally, $n\ge k$ and $n+\ell\ge2k+c$ imply
\[
 q^\ell\frac{q^{n-1}-1}{q^{k-1}-1}
 \ge q^{n+\ell-k}
 \ge q^{k+c}
 >q^c\qbinom k1.
\]
Thus \eqref{eq:singular-condition} holds, proving the last assertion.
\end{proof}

\begin{corollary}[The case $\ell=2c$]\label{cor:singular-balanced}
Let $k\ge2$, $c\ge1$, $\ell=2c$, and
$n\ge\max\{k,2k-c\}$.
If families $\mathcal A\subseteq\mathcal F_1$ and
$\mathcal B\subseteq\mathcal F_2$ satisfy
$A\cap B\ne\{0\}$ for every
$A\in\mathcal A$ and $B\in\mathcal B$, then
\[
 |\mathcal A||\mathcal B|
 \le q^{2c(k-1)}\qbinom{n-1}{k-1}\qbinom{2c}{c}.
\]
Equality holds if and only if there exists a point
$P\nsubseteq W$ such that
$ \mathcal A=(\mathcal F_1)_P$ and $\mathcal B=(\mathcal F_2)_P $. 
\end{corollary}
\begin{proof}
Since $n\ge\max\{k,2k-c\}$, we have
$n+\ell=n+2c\ge2k+c$.
The result follows from Theorem~\ref{thm:singular}.
\end{proof}

\subsection{Comparison with general cross-product criteria}\label{subsec:borg}
We compare Theorem~\ref{thm:singular} with two set-system criteria
under the natural outside-point encoding. Replace each
$U\in\mathcal F_a$ by its set of points outside $W$, obtaining
$\widehat{\mathcal F}_a$. These points span $U$, so the encoding
is injective. It also preserves the cross-intersection condition,
since every member of $\mathcal F_1$ is disjoint from $W$. 
The resulting set families are $r_a$-uniform, with
\[
 r_1=\qbinom k1,\qquad r_2=q^c.
\]
For a set family $\mathcal H$ on a finite ground set $\Omega$, define
\[
 \Delta_j(\mathcal H)=\max_{S\in\binom\Omega j}
       |\{H\in\mathcal H:S\subseteq H\}|,
\]
and set $\Delta_j(\mathcal H)=0$ for $j>|\Omega|$. Then
\begin{equation}\label{eq:encoded-ratio}
 \Delta_1(\widehat{\mathcal F}_2)=\qbinom\ell c,
 \qquad
 \Delta_2(\widehat{\mathcal F}_2)=\qbinom{\ell-1}{c-1}.
\end{equation}
Indeed, two distinct points outside $W$ can lie in a common
member of $\mathcal F_2$ only if their images in $V/W$ coincide.
In that case, their span meets $W$ in a point $R$, and the
members of $\mathcal F_2$ containing both points correspond
to the $c$-subspaces of $W$ containing $R$.
Thus, when $\ell=2c$,
\[
 \frac{\Delta_1(\widehat{\mathcal F}_2)}
      {\Delta_2(\widehat{\mathcal F}_2)}
 =\frac{q^{2c}-1}{q^c-1}=q^c+1.
\]

\emph{Star ratios.}
For cross $1$-intersection, Borg's sufficient condition
\cite[Theorem~1.2]{Borg2017} is
$\Delta_1\ge(r_1r_2+1)\Delta_2$ for each family.
When $\ell=2c$, this condition fails for $\widehat{\mathcal F}_2$,
since $q^c+1<r_1r_2+1$.

\emph{Weak spreadness.}
A weakly $(\sigma,1)$-spread family satisfies
$\Delta_{1+j}\le\sigma^{-j}\Delta_1$ for all $j\ge0$.
Thus \eqref{eq:encoded-ratio} implies $\sigma\le q^c+1$
when $\ell=2c$.
Set $h=\max\{r_1,r_2\}$.
After ordering the two uniformities, the product theorem
of Wen and Lv \cite[Theorem~1.3(ii)]{WLUnified2026} requires
\[
 \sigma\ge2e\max\left\{9,\frac{h}{1-1/h}\right\}.
\]
This cannot hold, since
$q^c+1\le h+1<2eh/(1-1/h)$.

Nevertheless, Theorem~\ref{thm:singular} applies at $\ell=2c$
whenever $n\ge\max\{k,2k-c\}$.
By \eqref{eq:singular-d2}, among the $r_1$ points of a fixed
$A\in\mathcal F_1$, at most one maps to the prescribed point of $V/W$.
The singular-space results in \cite{OLW2014,ZLL2016} concern
one-family EKR or semilattice problems and do not give this mixed
cross-product conclusion. In particular,
Corollary~\ref{cor:singular-balanced} is not a direct consequence of
\cite{OLW2014,ZLL2016}. It is also not covered, under the
outside-point encoding used here, by Borg's criterion
\cite[Theorem~1.2]{Borg2017} or the weak-spreadness criterion of
Wen and Lv \cite[Theorem~1.3(ii)]{WLUnified2026}, since the required
inequalities fail as shown above. This comparison does not address
other encodings or methods.

\subsection{A dual-polar example}\label{subsec:dual-polar}
Let $V$ be an eight-dimensional orthogonal space of minus type over $\Fq$.
The generators of $Q^-(7,q)$ are its three-dimensional totally singular
subspaces. Write
\[
 V=E\perp A,
\]
where $E$ is a six-dimensional hyperbolic space and $A$ is an anisotropic
plane. Let $X_r$ denote the set of $r$-dimensional totally singular
subspaces of $V$, and set
\[
 \mathcal F=\{G\in X_3:\dim(G\cap E)=1\},\qquad
 \mathcal P=\{S\in X_2:S\cap E=\{0\}\}.
\]

\begin{corollary}\label{cor:dual-polar-type}
Let $\mathcal A,\mathcal B\subseteq\mathcal F$ be families such that 
the intersection $G\cap H$ contains a member of $\mathcal P$ for all $G\in\mathcal A$ and $H\in\mathcal B$. Then
\[
 |\mathcal A||\mathcal B|\le(q^2+1)^2.
\]
Equality holds if and only if there exists $S\in\mathcal P$ such that
$\mathcal A=\mathcal B=\mathcal F_S$.
\end{corollary}

\begin{proof}
We apply Theorem~\ref{thm:criterion} with
$Y_1=Y_2=\mathcal F$. Every $G\in\mathcal F$ contains exactly $q^2$
members of $\mathcal P$, namely the two-dimensional complements of
$G\cap E$ in $G$. Thus $v_1=v_2=q^2$.

Fix $S\in\mathcal P$. Projection from $V=E\perp A$ onto $A$ restricts
to an isomorphism on $S$. The $E$-component of $S$ is therefore an
anisotropic plane $C$ anti-isometric to $A$. Writing $E=C\perp D$, the
space $D$ is a four-dimensional orthogonal space of minus type. The
members of $\mathcal F$ containing $S$ are precisely the spaces
$S+R$, where $R$ is a singular point of $D$. Hence
$L_1=L_2=q^2+1$.

It remains to bound the local count. Fix $S\in\mathcal P$ and
$Y\in\mathcal F$ with $S\nsubseteq Y$. If $S\cap Y=\{0\}$, then for every $T\in\mathcal F_S$,
we have $\dim(T\cap Y)\le 1$, since $S\subseteq T$
and $\dim T=3$. Hence the local count is zero.  Otherwise, put $R=S\cap Y$, which is a point.
For a counted $T\in\mathcal F_S$, the plane
$U=T\cap Y$ belongs to $\mathcal P$ and contains $R$.
Since $T=S+U$ is totally singular, we have
$U\le Y\cap S^\perp$.
If $\dim(Y\cap S^\perp)\ge3$, then
$Y\le S^\perp$, so $S+Y$ would be a four-dimensional
totally singular subspace, a contradiction.
Hence $\dim(Y\cap S^\perp)\le2$, and there is at most
one possible $U$. Thus we may take $d_1=d_2=1$.

Since $q^2<q^2+1$, the hypotheses of
Theorem~\ref{thm:criterion} hold, and the result follows.   
\end{proof}

\begin{remark}\label{rem:dual-polar-literature}
Ihringer \cite{Ihringer2015} treats cross-intersecting families in the
full set of generators of a classical polar space, with every cross pair
required to meet in at least a point. Ihringer and Metsch
\cite{IhringerMetsch2018} study large $\{0,1,\ldots,t\}$-cliques in dual
polar graphs. Neither result imposes the type restriction
$\dim(G\cap E)=1$ together with the permitted-centre condition
$S\cap E=\{0\}$. Thus Corollary~\ref{cor:dual-polar-type} is not a direct
consequence of the cited theorems. We make no broader priority claim.
\end{remark}

\subsection{Four-dimensional subspaces of different form ranks}\label{subsec:symplectic}
Let $V$ be a nondegenerate symplectic space of dimension $2\nu$ over
$\Fq$, with alternating form $f$. For $U\le V$, write
\[
 \rad U=\{u\in U:f(u,v)=0\text{ for all }v\in U\}.
\]
Let $\mathcal G_1$ and $\mathcal G_2$ be the families of four-dimensional
subspaces of $V$ on which $f$ has rank four and two, respectively.
The distinguished class $\mathcal P$ consists of the nondegenerate
planes. Thus members of the two families have the same dimension
but different form ranks. This differs from the totally isotropic
setting of \cite{YWPolar}.
The enumeration of subspaces in symplectic spaces is studied in
\cite{GLW2013}. We derive the required counts and radical-dependent
local estimates below.

\begin{theorem}\label{thm:mixed-symplectic}
Suppose $\nu\ge5$.
Let $\mathcal A\subseteq\mathcal G_1$ and
$\mathcal B\subseteq\mathcal G_2$ be families such that
$A\cap B$ contains a nondegenerate plane for every
$A\in\mathcal A$ and $B\in\mathcal B$. Then
\begin{equation}\label{eq:mixed-symp-bound}
 |\mathcal A||\mathcal B|
 \le q^{2\nu-4}\qbinom{2\nu-4}1
       \left(\frac{q^{2\nu-2}-1}{q^2-1}\right)^2.
\end{equation}
Equality holds if and only if there exists a nondegenerate
plane $P$ such that
$\mathcal A=(\mathcal G_1)_P$ and $\mathcal B=(\mathcal G_2)_P$.
\end{theorem}

\begin{proof}
Put
\[
 \rho=\frac{q^{2\nu-2}-1}{q^2-1},\qquad
 \eta_1=q^{2\nu-4},\qquad \eta_2=\qbinom{2\nu-4}1.
\]
We apply Theorem~\ref{thm:criterion} with $Y_a=\mathcal G_a$
for $a=1,2$ and the following parameters. 
\begin{equation}\label{eq:symp-counts}
 \renewcommand{\arraystretch}{1.25}
 \begin{array}{c|cc}
  &\mathcal G_1&\mathcal G_2\\ \hline
 v_a&q^2(q^2+1)&q^4\\
 L_a&\eta_1\rho&\eta_2\rho\\
 d_a\text{ (upper bound)}&q^2\eta_1&q^2\eta_2
 \end{array}
\end{equation}
We justify the incidence numbers, star sizes, and local bounds.

\emph{Incidence numbers.}
In a nondegenerate $d$-dimensional symplectic space with $d\ge4$,
there are $(q^d-1)q^{d-1}$ ordered pairs $(u,v)$ satisfying
$f(u,v)=1$.
Each nondegenerate plane contains $(q^2-1)q$ such pairs,
so the number of nondegenerate planes is
$q^{d-2}(q^d-1)/(q^2-1)$.
Counting incidences between points and totally isotropic planes
gives $\qbinom{d-2}1(q^d-1)/(q^2-1)$ totally isotropic planes:
each point lies in $\qbinom{d-2}1$ such planes, and each plane
contains $q+1$ points.

The first formula with $d=4$ gives $v_1$.
A four-space of form rank two is the direct sum of a nondegenerate
plane and its two-dimensional radical.
Its nondegenerate planes are precisely the complements to the
radical, giving $v_2=q^4$.
Every four-space containing $P\in\mathcal P$ has a unique
decomposition $P\perp Q$, where $Q\le P^\perp$ is a plane.
The two formulas with $d=2\nu-2$ give $L_1$ and $L_2$,
respectively.

Any three-space $S$ containing $P$ has the form $P\perp R$,
where $R\le P^\perp$ is a point.
The members of $\mathcal G_1$ and $\mathcal G_2$ containing $S$
correspond to the nondegenerate and totally isotropic planes
of $P^\perp$ through $R$, respectively.
Their numbers are $\eta_1$ and $\eta_2$.
For the first count, subtract the number of points in $R^\perp/R$
from that in $P^\perp/R$; for the second, count the points
of $R^\perp/R$.
Here $R^\perp$ is taken within $P^\perp$.

\emph{Local counts.}
Fix $a\in\{1,2\}$, $P\in\mathcal P$, and
$Y\in\mathcal G_{\bar a}$ with $P\nsubseteq Y$.
We count the members $T\in\mathcal G_a$ containing $P$
for which $D=T\cap Y$ contains a nondegenerate plane.
If $P\cap Y=\{0\}$, any nondegenerate plane $Q\le D$
determines $T=P+Q$.
Thus there are at most $v_{\bar a}$ such $T$, and
$v_{\bar a}\le q^2\eta_a$ for $\nu\ge5$.

Suppose now that $R=P\cap Y$ is a point.
Since $T$ and $Y$ have different form ranks, $\dim D\le3$.
If $D\le E$ and $E$ has form rank two, then any nondegenerate
plane $Q_0\le D$ gives
\[
 E=Q_0\perp\rad E,\qquad
 D=Q_0\perp(D\cap\rad E),\qquad
 \rad D=D\cap\rad E.
\]
If $a=2$, apply this with $E=T$.
Since $R\le P$ and $P\cap\rad T=\{0\}$,
we have $R\nsubseteq\rad D$.
If $a=1$, apply it with $E=Y$.
Then $R\le\rad D$ if and only if $R\le\rad Y$.

Suppose first that $R\nsubseteq\rad D$.
Choose a nonzero $r\in R$ and $v\in D$ with $f(r,v)\ne0$.
Then $Q=\langle r,v\rangle$ is a nondegenerate plane
in $D$ containing $R$.
Since $R\nsubseteq\rad Y$, the number of nondegenerate
planes of $Y$ containing $R$ is
$\qbinom31-\qbinom21=q^2$.
For each such $Q$, the three-space $P+Q$ is contained
in exactly $\eta_a$ members of $\mathcal G_a$.
Thus there are at most $q^2\eta_a$ possible $T$.

It remains to consider $a=1$ with $R\le\rad Y$.
Then $\dim D=3$ and $D\cap\rad Y=R$.
Consequently, $D/R$ is a complement to the point
$(\rad Y)/R$ in $Y/R$.
There are $q^2$ such complements, each of which determines
$D$ and hence $T=P+D$.
This gives at most $q^2\le q^2\eta_1$ possible $T$.
The local bounds in \eqref{eq:symp-counts} follow.

Finally, $\nu\ge5$ implies $\rho\ge1+q^2+q^4+q^6$, so
\[
 v_1d_1\le(q^6+q^4)\eta_1<\rho\eta_1=L_1,
 \qquad
 v_2d_2\le q^6\eta_2<\rho\eta_2=L_2.
\]
The conclusion follows from Theorem~\ref{thm:criterion}.
\end{proof}

\begin{corollary}[The case $\nu=5$]\label{cor:mixed-symplectic-five}
Let $V$ be a nondegenerate symplectic space of dimension $10$ over
$\Fq$, and let $\mathcal G_1$ and $\mathcal G_2$ be as above.
If families $\mathcal A\subseteq\mathcal G_1$ and
$\mathcal B\subseteq\mathcal G_2$ satisfy
$A\cap B$ contains a nondegenerate plane for every
$A\in\mathcal A$ and $B\in\mathcal B$, then
\[
 |\mathcal A||\mathcal B|
 \le q^6\qbinom61
 \left(\frac{q^8-1}{q^2-1}\right)^2.
\]
Equality holds if and only if there exists a nondegenerate plane $P$
such that
$\mathcal A=(\mathcal G_1)_P$ and $\mathcal B=(\mathcal G_2)_P$.
\end{corollary}
\begin{proof}
This is Theorem~\ref{thm:mixed-symplectic} with $\nu=5$.
\end{proof}

\begin{remark}\label{rem:symp-comparison}
At $\nu=5$, the condition in Corollary~\ref{cor:balanced} fails:
$\rho=1+q^2+q^4+q^6<v_1v_2=q^8+q^6$, and the maximum number of
members of $\mathcal G_a$ containing a three-space is at least $\eta_a$.
The direct local estimates nevertheless satisfy the hypotheses of
Theorem~\ref{thm:criterion}.
For $\nu\ge6$, Borg's star-ratio criterion
\cite[Theorem~1.2]{Borg2017}, applied to the encoding by nondegenerate
planes, also gives the product bound, since its required inequality becomes
$\rho\ge v_1v_2+1$. At $\nu=5$, this inequality fails. Thus
Corollary~\ref{cor:mixed-symplectic-five} is not covered by Borg's
criterion under this encoding. The affine and symplectic-polar
results \cite{YWAffine,YWPolar} concern different geometric settings and
do not imply the mixed form-rank $(4,2)$ result here. We did not locate the
same $\nu=5$ conclusion in the references surveyed for this paper.
Moreover, the unrestricted Grassmann $2$-star through a nondegenerate
plane has size $L_1+L_2$, so ignoring the type restrictions gives
the larger star product $(L_1+L_2)^2$ rather than $L_1L_2$.
\end{remark}

\begin{corollary}[An equal-type three-space case]\label{cor:symp-three}
Let $\nu\ge3$, and let $\mathcal T$ consist of the three-dimensional
subspaces of a symplectic $2\nu$-space on which the form has rank two.
If $Z_1,Z_2\subseteq\mathcal T$ are families such that
$T_1\cap T_2$ contains a nondegenerate plane for every
$T_1\in Z_1$ and $T_2\in Z_2$, then
\[
 |Z_1||Z_2|\le\qbinom{2\nu-2}1^{\!2}.
\]
Equality holds if and only if there exists a nondegenerate plane
$P$ such that $Z_1=Z_2=\mathcal T_P$.
\end{corollary}
\begin{proof}
Here $v=q^2$, $L=\qbinom{2\nu-2}1$, and the same local
count as above gives $d=q$. Since $\nu\ge3$, we have
$vd=q^3<L$, and Theorem~\ref{thm:criterion} applies.  
\end{proof}

\section{Standard full-fibre consequences}\label{sec:applications}
Table~\ref{tab:standard-examples} records the coarse extension criterion
for standard semilattices \cite{Delsarte1976,QR1997,Suda2012}.
Assume $1\le s<N_1\le N_2$ and put
\[
 v_a=\nu(s,N_a),\qquad L_a=\theta(s,N_a),\qquad
 R_a=\frac{\theta(s,N_a)}{\theta(s+1,N_a)}
 \qquad(a=1,2).
\]
Suppose that
\begin{equation}\label{eq:mixed-ratio}
 R_1>v_1v_2\qquad\text{and}\qquad R_2>v_1v_2.
\end{equation}
Then Corollary~\ref{cor:balanced} gives
\[
 |Z_1||Z_2|\le L_1L_2
\]
for all cross $s$-intersecting families
$Z_a\subseteq X_{N_a}$.
Equality holds if and only if there exists $x\in X_{s}$ such that $Z_a=(X_{N_a})_x$  for $a=1,2$.
In every row, $R_1\ge R_2$, so it is enough to check
$R_2>v_1v_2$.
Write
\[
 C=\binom{N_1}s\binom{N_2}s,\qquad
 C_q=\qbinom{N_1}s\qbinom{N_2}s,\qquad
 \Delta=N_1+N_2-2s.
\]

\begin{table}[htbp]
\centering
\small
\setlength{\tabcolsep}{4pt}
\renewcommand{\arraystretch}{1.5}
\begin{tabular}{@{}lcc@{}}
\toprule
Semilattice & Product bound $B=L_1L_2$ & Sufficient condition \\
\midrule
Boolean (Johnson)
 & $\displaystyle\prod_{a=1}^2\binom{n-s}{N_a-s}$
 & $\displaystyle\frac{n-s}{N_2-s}>C$ \\[4pt]
Grassmann
 & $\displaystyle\prod_{a=1}^2\qbinom{n-s}{N_a-s}$
 & $\displaystyle\frac{q^{n-s}-1}{q^{N_2-s}-1}>C_q$ \\[4pt]
\shortstack[l]{Hamming\\(partial words)}
 & $\displaystyle h^{\Delta}\prod_{a=1}^2\binom{m-s}{N_a-s}$
 & $\displaystyle\frac{h(m-s)}{N_2-s}>C$ \\[4pt]
\shortstack[l]{Bilinear forms\\(partial linear maps)}
 & $\displaystyle q^{\beta\Delta}\prod_{a=1}^2\qbinom{m-s}{N_a-s}$
 & $\displaystyle q^\beta\frac{q^{m-s}-1}{q^{N_2-s}-1}>C_q$ \\[4pt]
Partial injections
 & $\displaystyle\prod_{a=1}^2\left(\binom{m-s}{N_a-s}
        \frac{(n-s)!}{(n-N_a)!}\right)$
 & $\displaystyle\frac{(m-s)(n-s)}{N_2-s}>C$ \\[4pt]
Signed sets
 & $\displaystyle h^{\Delta}\prod_{a=1}^2\binom{n-s}{N_a-s}$
 & $\displaystyle\frac{h(n-s)}{N_2-s}>C$ \\[4pt]
\shortstack[l]{Restricted\\signed sets}
 & $\displaystyle(n-1)^{\Delta}\prod_{a=1}^2\binom{n-s}{N_a-s}$
 & $\displaystyle\frac{(n-1)(n-s)}{N_2-s}>C$ \\
\bottomrule
\end{tabular}
\smallskip
\caption{Mixed-rank bounds $|Z_1||Z_2|\le B$ for $1\le s<N_1\le N_2$.
Each sufficient condition is $R_2>v_1v_2$; equality means full stars at a
common centre. Here $v_1v_2=C_q$ in the Grassmann and bilinear-form rows,
and $C$ otherwise.}
\label{tab:standard-examples}
\end{table}

Here $h\ge1$ and $\beta\ge0$ are integers, and
$[d]=\{1,\ldots,d\}$. In the map examples, the order is restriction and
the rank is the size or dimension of the domain; cross $s$-intersection
means that each pair from the two families has a common restriction of
rank at least $s$. The Hamming and signed-set rows describe the same
construction, while the restricted signed-set row corresponds to
$h=n-1$ after relabelling the allowed signs. See \cite{Suda2012} for
details of the standard semilattices and their incidence parameters.

These examples are included as applications of the criterion, not as
separate novelty claims. Cross-intersection results are already known
for signed sets and partial permutations \cite{BL2010,Borg2010}, and
stronger bounds are known for sets and Grassmann spaces
\cite{BaoJi2026,CLLW2023,FLST2014,ST2014}.
The bilinear-forms row of Table 1 also covers attenuated spaces. Indeed, let \(V=V_0\oplus W\), where \(\dim V_0=n\) and \(\dim W=\ell\). Every subspace \(U\le V\) satisfying \(U\cap W=\{0\}\) is the graph of a linear map from a subspace of \(V_0\) to \(W\). Thus the corresponding cross-intersection result follows from Table 1 by taking \(m=n\) and \(\beta=\ell\).

\section*{Acknowledgements}
The author acknowledges the support of JSPS KAKENHI Grant Number 26K06904.

\section*{AI Declaration}
OpenAI ChatGPT (GPT-5.6 Sol and GPT-6) was used for literature searches, English-language editing,
reorganization, and the development and checking of arguments, including
the mixed-type incidence estimates. The author is responsible for the
mathematical content, verification of references, and final submission.

\end{document}